\input amstex
\documentstyle{amsppt}
\magnification 1200
\vcorrection{-1cm}
\input epsf
\NoBlackBoxes

\rightheadtext{ Separating semigroup of hyperelliptic and genus 3 curves }
\topmatter
\title
         Separating semigroup of hyperelliptic curves and of genus 3 curves
\endtitle

\author  S.~Yu.~Orevkov
\endauthor

\address
Steklov Mathematical Institute, Gubkina 8, Moscow, Russia
\endaddress

\address
IMT, l'universit\'e Paul Sabatier, 118 route de Narbonne, Toulouse, France
\endaddress

\email
orevkov\@math.ups-tlse.fr
\endemail

\endtopmatter

\def\NN{\Bbb N}
\def\R {\Bbb R}

\def\P {\Bbb P}
\def\RP {\Bbb{RP}}
\def\CP {\Bbb{CP}}
\def\sep{\operatorname{Sep}} \def\Sep{\sep}
\def\ch{\operatorname{ch}}
\def\conj{\operatorname{conj}}
\def\sign{\operatorname{sign}}
\def\Sym{\operatorname{Sym}}

\def\refA   {1}
\def\refKS  {2}
\def\refN   {3}

\def\thQ   {1}
\def\thHE  {2}

\def\sectVDM {2}
\def\propVDM  {\sectVDM.1}
\def\corVDM   {\sectVDM.2}

\def\sectHE {3}
\def\lemSpe {\sectHE.1}
\def\lemAJ  {\sectHE.2}

\def\eqVDM   {1}
\def\eqAJ    {2}

\document

\head 1. Introduction
\endhead

By a real algebraic curve we mean a complex algebraic curve $C$
endowed with an antiholomorphic involution $\conj:C\to C$
(the complex conjugation involution). In this case
we denote the {\it real locus} $\{p\in C\mid \conj(p)=p\}$ by $\R C$.
A real curve is of {\it dividing type} (or of {\it type I\/}) if $\R C$
divides $C$ into two halves excanged by the complex conjugation.
All curves considered here are smooth and irreducible.

A sufficient condition for $C$ to be of dividing type is the existence
of a {\it separating morphism} $f:C\to\P^1$, that is a morphism such that
$f^{-1}(\RP^1)=\R C$.
It follows from Ahlfors' results [\refA] that this condition is also necessary:
any real curve of dividing type admits a separating morphism.
The restriction of a separating morphism to $\R C$ is a covering over $\RP^1$.
If we fix the numbering of connected components $A_1,\dots,A_n$ of $\R C$,
we may consider the sequence $d(f)=(d_1,\dots,d_n)$ where $d_i$ is
the covering degree of $f$ restricted to $A_i$.

Kummer and Shaw studied in [\refKS] the following problem. Given a curve $C$
of dividing type, which sequences are realizable as $d(f)$ for separating morphisms
$f:C\to\P^1$? It is easy to see that the set of all realizable sequences is an
additive semigroup (see [\refKS, Proposition 2.1]). Following [\refKS], we call it
the {\it separating semigroupe} of $C$ and denote by $\Sep(C)$.

Several interesting properties of $\Sep(C)$ are established in [\refKS].
In particular, it is shown that $\Sep(C)=\NN^{g+1}$ for an $M$-curve $C$
(a curve $C$ of genus $g$ is called an $M$-curve
if $\R C$ has $g+1$ connected components which is the maximal possible number for genus
$g$ curves).
Also it is shown in [\refKS] that sometimes the separating semigroup does depend on
the numbering of the components. The simplest example is a hyperbolic quartic curve in
$\RP^2$ (a plane curve is called hyperbolic if the linear projection from some point is
a separating morphism).
If $C$ is such a curve, then $\R C$ consists of two ovals one inside another.
If we number them so that the inner oval is first, then we have $(1,2)\in\Sep(C)$
but $(2,1)\not\in\Sep(C)$; see [\refKS, Example 3.7]. Moreover, $\Sep(C)$ is
almost computed in [\refKS]: it is shown that $\NN\times\NN_{\ge2}\subset\Sep(C)$.
We complete this computation:

\proclaim{ Theorem \thQ }
Let $C$ be a nonsingular real hyperbolic quartic curve in $\RP^2$ whose
ovals are numbered so that the inner one is first.
Then $\sep(C) = \NN\times\NN_{\ge2}$.
\endproclaim

The proof relies on two facts: Theorem \thHE\ below and Natanzon's theorem [\refN, Theorem 2.3]
which states that two branched coverings over a disk are (left-right) topologically
equivalent if and only if they are equivalent over the boundary circle.

\proclaim{ Theorem \thHE }
Let $C$ be a real hyperelliptic curve of genus $g\ge 2$ of dividing type
but not an $M$-curve.
Then
$$
   \sep(C) = \cases
          (1,1)\NN \cup (\NN_{\ge(g+1)/2})^2
                   &\text{if $g$ is odd,}\\
         2\NN \cup \NN_{\ge g}
     &\text{if $g$ is even.}
      \endcases
$$
\endproclaim

Note that any real genus three curve of dividing type is either an $M$-curve,
or hyperelliptic, or a plane hyperbolic quartic. Thus the results of [\refKS]
completed by our Theorems 1 and 2 provide separating semigroups of all
real curves of dividing type up to genus $3$.

\head\sectVDM. Dual Vandermonde system of equations
\endhead

Let $x_1,x_2,\dots,x_n$ be real numbers.
We consider the homogeneous system of linear equations
with indeterminates $h_1,\dots,h_n$ (the dual Vandermonde system):
$$
     \sum_{i=1}^n x_i^k h_i = 0, \qquad  k=0,\dots,g-1.    \eqno(\eqVDM)
$$
This condition on $(h_1,\dots,h_n)$ can be equivalently rewritten as follows
$$
     \sum_{i=1}^n h_i F(x_i) = 0 \qquad\text{for any $F\in\R[x]$ with $\deg F<g$}.
$$
Given a sequence of real numbers $h=(h_1,\dots,h_n)$, we define $\ch(h)$ as
{\it the number of changes of sign} of $h$, i.~e., the number of pairs $(i,j)$ such
that $1\le i<j\le n$, $h_i h_j<0$, and $h_k=0$ if $i<k<j$.

\proclaim{ Proposition \propVDM } Let $x_1<\dots<x_n$, $n>0$.
A sequence $s=(s_1,\dots,s_n)$ with $s_i\in\{-1,0,1\}$ is
the sequence of signs of a non-zero solution to the system (\eqVDM)
if and only if $\ch(s)\ge g$.
\endproclaim

\demo{ Proof }
($\Rightarrow$). Suppose that $h=(h_1,\dots,h_n)$ is a solution to (\eqVDM), and $\ch(h)<g$.
Then we can choose a polynomial $F$ of degree less than $g$ such that
$F(x_i)\ne 0$ and $h_iF(x_i) \ge 0$ for any $i=1,\dots,n$.
Then $h_1F(x_1)+\dots+h_n F(x_n)=0$ and each term in this sum is non-negative.
Hence $h=(0,\dots,0)$.

\smallskip
($\Leftarrow$). Let $\ch(s)\ge g$.
Let $I=\{i_0,\dots,i_g\}\subset\{1,\dots,n\}$ be
such that $\ch(s_{i_0},\dots,s_{i_g})\break=g$.
Let $h'_I=(h'_{i_0},\dots,h'_{i_g})$ be a non-zero solution to the system (\eqVDM)
with $\sum_{1\le i\le n}$ replaced by $\sum_{i\in I}$.
By ``($\Rightarrow$)'' part, we have $\ch(h'_I)=g$.
Thus, changing the sign of $h'_I$ if necessary,
we have $\sign h'_{i_j}=s_{i_j}\ne 0$ for all $j=0,\dots,g$.

Let us choose $(h_i)_{i\not\in I}$ such that $\sign h_i=s_i$ and $|h_i|<\varepsilon\ll 1$.
Set $h_{i_0}=h'_{i_0}$. Since the Vandermonde $(g\times g)$-determinant
(corresponding to the columns numbered by $I\setminus\{i_0\}$)
is non-zero, the remaining numbers $h_{i_1},\dots,h_{i_g}$ are uniquely determined by (\eqVDM).
Moreover, if $\varepsilon$ is small enough, then $h_I=(h_i)_{i\in I}$ is close to $h'_I$,
thus $\sign h_i=\sign h'_i=s_i$ for all $i\in I$.
\qed\enddemo

\proclaim{ Corollary \corVDM }
Let $x_1,\dots,x_n$ be real numbers, not necessarily distinct. For $x\in\R$
we set $I(x)=\{i\mid x_i=x\}$. Let
$(h_1,\dots,h_n)$ be a real solution to the system (\eqVDM) such that
$h_i\ne 0$ for all $i=1,\dots,n$. Then at least one of the following
two possibilities takes place:
\smallskip
\roster
\item"(i)" $\sum_{i\in I(x)}h_i=0$ for any $x\in\R$, in particular, each $x_i$
           occurs at least twice in the sequence $(x_1,\dots,x_n)$;
\smallskip
\item"(ii)" the sequence $(h_1,\dots,h_n)$ contains at least $[(g+1)/2]$ positive
and at least $[(g+1)/2]$ negative members.
\endroster
\endproclaim

\head 3. Separating semigroup of hyperelliptic curves
\endhead
In this section we prove Theorem \thHE.

\proclaim{ Lemma \lemSpe }
Let $C$ be a (complex) hyperelliptic curve of genus $g$
and $f$ a meromorphic function on $C$ such that the zero divisor $(f)_0$
is special (this is so, for example, when $\deg f<g$). Then
$f=f_1\circ\pi$ where $\pi:C\to\P^1$ is the hyperelliptic projection
and $f_1$ a meromorphic function on $\P^1$.
\endproclaim

\demo{ Proof }
If $D$ and $D'$ are two effective divisors on a curve,
then the embedding $\phi_D$ defined by the complete linear system $|D|$ is a
composition of $\phi_{D+D'}$ with a linear projection.
Let $D=(f)_0$ and let $D'$ be an effective divisor such that $D+D'\sim K_C$ (such $D'$ exists
since $D$ is special). Thus $\phi_D$ is a projection of
the canonical embedding which is known to factor through the hyperelliptic
projection.
\qed\enddemo

\proclaim{ Lemma \lemAJ }
Let $C$ be a real algebraic curve of genus $g>0$ of dividing type.
Let $\omega_1,\dots,\omega_g$ be a base of holomorphic
1-forms on $C$.

\smallskip
(a).
Let  $f:C\to\P^1$ be a separating morphism and $\{p_1,\dots,p_n\}=f^{-1}(p)$
for some point $p\in\RP^1$. Then there exist real positive
(with respect to a fixed complex orientation) tangent vectors $v_1,\dots,v_n$
($v_i$ tangent at $p_i$) such that
$$
   \sum_{i=1}^n \omega_k(v_i)=0 \qquad\text{for each $k=1,\dots,g$}.     \eqno(\eqAJ)
$$

(b). Conversely, let $p_1,\dots,p_n$ be distinct points on $\R C$ and
$v_1,\dots,v_n$ be positive real tangent vectors ($v_i$ tangent at $p_i$) such that
(\eqAJ) holds.
Suppose in addition that the divisor $D=p_1+\dots+p_n$ is non-special
 i.~e., $h^0(K_C-D)=0$. Then there exists
a separating morphism with fiber $D$.
\endproclaim

\demo{ Proof } (a). Follows from Abel-Jacobi Theorem.

\smallskip
(b). Follows from Abel-Jacobi Theorem combined with [\refKS, Lemma 2.10].
Indeed, consider the Abel-Jacobi mapping $\varphi:\Sym^n(C)\to\Cal J(C)$.
The condition (\eqAJ) means that $v=(v_1,\dots,v_n)$ considered as a tangent
vector to $\Sym^n(C)$ at $D$ is in the kernel of the differencial of $\varphi$ at $D$.
The non-specialness of $D$ means that $\varphi$ is a submersion near $D$, hence
$v$ is tangent to $\varphi^{-1}(\varphi(D))=|D|$ at $D$.
Hence there exists a path $[0,t_0]\to|D|$, $t\mapsto D_t$, such that $D_0=D$ and
$\big({d\over dt}D_t\big)_{t=0}=v$. Then, for any $t$, $0<t\le t_0$, there exists a
meromorphic function $f_t:C\to\P^1$ such that $D=(f_t)_0$ and $D_t=(f_t)_\infty$.
If $t$ is small enough, then the condition of positivity of the $v_i$'s implies that
the zeros and poles of $f_t$ interlace along $\R C$, thus
$f_t$ is a separating morphism by [\refKS, Lemma 2.10].
\qed\enddemo

\demo{ Proof of Theorem \thHE}
Let $C$ be a real hyperelliptic curve of genus $g\ge 2$ which is not an $M$-curve.
Then it is given by an equation $y^2=G(x)$ where $G(x)$ is a generic real polynomial
of degree $2g+2$ positive everywhere on $\R$.
We consider the standard base of holomorphic 1-forms
$\omega_1,\dots,\omega_g$ where $\omega_k=x^{k-1}dx/y$.
The hyperelliptic projection is given by $(x,y)\mapsto x$.
Its restriction to $\R C$ is an unramified two-fold covering over $\RP^1$
which is trivial for even $g$ and non-trivial for odd $g$.
We choose the complex orientation on $\R C$ such that $dx>0$ on positive
tangent vectors.

Let $f:C\to\P^1$ be a separating morphism and let
$\{p_1,\dots,p_n\}=f^{-1}(p)$ for a generic $p\in\RP^1$. We set $p_i=(x_i,y_i)$, $i=1,\dots,n$.
By Lemma \lemAJ(a) there exist positive tangent vectors $v_1,\dots,v_n$
such that (\eqAJ) holds. Let $a_i=dx(v_i)$. The positivity of $v_i$ means $a_i>0$.
Then (\eqAJ) takes the form (\eqVDM) for $h_i=a_i/y_i$, and
Theorem \thHE\ follows from Corollary \corVDM\ and  Lemma \lemSpe.
\enddemo

\head 4. Separating semigroup of genus three curves
\endhead

In this section we prove Theorem \thQ.
Let $C$ be a plane hyperbolic quartic curve. We have $\NN\times\NN_{\ge2}\subset\Sep(C)$,
see [\refKS, Example 3.7]. Let us prove the inverse inclusion.

It is shown in [\refKS, Example 2.8] thet $(1,1)\not\in\Sep(C)$.
Suppose there exists a separating morphism $f_0:C\to\P^1$ with $d(f_0)=(n,1)$, $n\ge 2$.
Let $C^+$ be one of the two halves into which $\R C$ divides $C$. Then the restriction
of $f_0$ to $C^+$ is a branched covering over a disk $\Delta$ which is one of the halves of
$\CP^1\setminus\RP^1$.
By perturbing $f_0$ (together with $C$) we may assume
that all critical values are simple, i.~e., $f^{-1}(p)$ has at least $n$ points
for any $p\in\Delta$.

Let $f_1:C\to\P^1$ be a separating morphism with $d(f_1)=(1,n)$ which exists by
[\refKS, Example 3.7]. It can be chosen so that all its critical values are simple.
Then, by Natanzon's result [\refN, Theorem 2.3], there exists a continuous family of branched
coverings $f_t:C^+\to\Delta$, $0\le t\le 1$, which connects $f_0$ with $f_1$.
Let $C^+_t$ be $C^+$ endowed with the complex structure lifted from $\Delta$ by $f_t$,
and let $C_t$ be $C^+_t$ glued along the boundary with its complex conjugate copy.
Then $f_t$ extends to a separating morphism $C_t\to\P^1$ which we also denote by $f_t$.
So, we obtain a continuous family of separating morphisms $f_t$ of genus three
curves $C_t$.

By continuity, we have $d(f_t)=(1,n)$ for a suitable numbering of the components of $\R C_t$.
Hence, by Theorem \thHE, the curve $C_t$ cannot be hyperelliptic for any $t$.
It is well-known that any non-hyperelliptic genus three curve is isomorphic
to a smooth quartic curve in $\P^2$. Thus there exists a continuous family
of embeddings $\iota_t:C_t\to\P^2$ such that $\iota_t(C_t)$ is a smooth real quartic
curve, and we have a continuous family of separating morphisms 
of them onto $\P^1$. The interior and exterior ovals cannot interchange in this family
which contradicts the fact that $d(f_0)\ne d(f_1)$ and the embedding to $\P^2$ is unique up
to projective equivalence.

\enddocument